\documentclass[11pt,twoside]{article}
\usepackage{amsmath, amsthm, amscd, amsfonts, amssymb, graphicx, color}

\date{}

\begin{document}

\centerline{}

\centerline {\Large{\bf Construction of continuous controlled $K$-$g$-fusion frames }}
\centerline {\Large{\bf  in Hilbert spaces }}

\newcommand{\mvec}[1]{\mbox{\bfseries\itshape #1}}
\centerline{}
\centerline{\textbf{Prasenjit Ghosh}}
\centerline{Department of Pure Mathematics, University of Calcutta,}
\centerline{35, Ballygunge Circular Road, Kolkata, 700019, West Bengal, India}
\centerline{e-mail: prasenjitpuremath@gmail.com}
\centerline{}
\centerline{\textbf{T. K. Samanta}}
\centerline{Department of Mathematics, Uluberia College,}
\centerline{Uluberia, Howrah, 711315,  West Bengal, India}
\centerline{e-mail: mumpu$_{-}$tapas5@yahoo.co.in}

\newtheorem{Theorem}{\quad Theorem}[section]

\newtheorem{definition}[Theorem]{\quad Definition}

\newtheorem{theorem}[Theorem]{\quad Theorem}

\newtheorem{remark}[Theorem]{\quad Remark}

\newtheorem{corollary}[Theorem]{\quad Corollary}

\newtheorem{note}[Theorem]{\quad Note}

\newtheorem{lemma}[Theorem]{\quad Lemma}

\newtheorem{example}[Theorem]{\quad Example}

\newtheorem{result}[Theorem]{\quad Result}
\newtheorem{conclusion}[Theorem]{\quad Conclusion}

\newtheorem{proposition}[Theorem]{\quad Proposition}

\begin{abstract}
\textbf{\emph{We present the notion of continuous controlled $K$-$g$-fusion frame in Hilbert space which is the generalization of discrete controlled $K$-$g$-fusion frame.\,We discuss some characterizations of continuous controlled $K$-$g$-fusion frame.\,Relationship between continuous controlled $K$-$g$-fusion frame and quotient operator is being studied.\,Finally, stability of continuous controlled $g$-fusion frame has been described.}}
\end{abstract}
\\
\\

\section{Introduction}
 
\smallskip\hspace{.6 cm}In 1952, Duffin and Schaeffer \cite{Duffin} introduced frame for Hilbert space to study some fundamental problems in non-harmonic Fourier series.\,Later on, after some decades, frame theory was popularized by Daubechies et al.\,\cite{Daubechies}.

Frame for Hilbert space was defined as a sequence of basis-like elements in Hilbert space.\,A sequence \,$\left\{\,f_{\,i}\,\right\}_{i \,=\, 1}^{\infty} \,\subseteq\, H$\, is called a frame for a separable Hilbert space \,$\left(\,H,\,\left<\,\cdot,\,\cdot\,\right>\,\right)$, if there exist positive constants \,$0 \,<\, A \,\leq\, B \,<\, \infty$\, such that
\[ A\; \|\,f\,\|^{\,2} \,\leq\, \sum\limits_{i \,=\, 1}^{\infty}\, \left|\ \left <\,f \,,\, f_{\,i} \, \right >\,\right|^{\,2} \,\leq\, B \,\|\,f\,\|^{\,2}\; \;\text{for all}\; \;f \,\in\, H.\]
For the past few years many other types of frames were proposed such as \,$K$-frame \cite{L}, fusion frame \cite{Kutyniok}, \,$g$-frame \cite{Sun}, \,$g$-fusion frame \cite{G, Ahmadi} and \,$K$-$g$-fusion frame \cite{Sadri} etc.\,P. Ghosh and T. K. Samanta \cite{Ghosh} have discussed generalized atomic subspaces for operators in Hilbert spaces.

Controlled frame is one of the newest generalization of frame.\,P. Balaz et al.\,\cite{B} introduced controlled frame to improve the numerical efficiency of interactive algorithms for inverting the frame operator.\,In recent times, several generalizations of controlled frame namely, controlled\,$K$-frame \cite{N}, controlled\,$g$-frame \cite{F}, controlled fusion frame \cite{AK}, controlled $g$-fusion frame \cite{HS}, controlled $K$-$g$-fusion frame \cite{GR} etc. have been appeared.\,Continuous frames were proposed by Kaiser \cite{Ka} and it was independently studied by Ali et al.\,\cite{Al}.\,At present, frame theory has been widely used in signal and image processing, filter bank theory, coding and communications, system modeling and so on.

In this paper,\, continuous controlled $K$-$g$-fusion frame in Hilbert spaces is studied and some of their properties are going to be established.\,Under some sufficient conditions, we will see that any continuous controlled $K$-$g$-fusion frame is equivalent to continuous $K$-$g$-fusion frame.\,A necessary and sufficient condition for continuous controlled $g$-fusion Bessel family to be a continuous controlled \,$K$-$g$-fusion frame with the help of quotient operator is established.\,At the end, we study some stability results of continuous controlled $g$-fusion frame.  

Throughout this paper,\;$H$\; is considered to be a separable Hilbert space with associated inner product \,$\left <\,\cdot \,,\, \cdot\,\right>$\, and \,$\mathbb{H}$\, is the collection of all closed subspace of \,$H$.\,$I_{H}$\; is the identity operator on \,$H$.\,$\mathcal{B}\,(\,H_{\,1},\, H_{\,2}\,)$\; is a collection of all bounded linear operators from \,$H_{\,1} \,\text{to}\, H_{\,2}$.\,In particular \,$\mathcal{B}\,(\,H\,)$\, denotes the space of all bounded linear operators on \,$H$.\;For \,$S \,\in\, \mathcal{B}\,(\,H\,)$, we denote \,$\mathcal{N}\,(\,S\,)$\; and \,$\mathcal{R}\,(\,S\,)$\, for null space and range of \,$S$, respectively.\,Also, \,$P_{M} \,\in\, \mathcal{B}\,(\,H\,)$\; is the orthonormal projection onto a closed subspace \,$M \,\subset\, H$.\,$\mathcal{G}\,\mathcal{B}\,(\,H\,)$\, denotes the set of all bounded linear operators which have bounded inverse.\,If \,$S,\, R \,\in\, \mathcal{G}\,\mathcal{B}\,(\,H\,)$, then \,$R^{\,\ast},\, R^{\,-\, 1}$\, and \,$S\,R$\, also belongs to \,$\mathcal{G}\,\mathcal{B}\,(\,H\,)$.\,$\mathcal{G}\,\mathcal{B}^{\,+}\,(\,H\,)$\, is the set of all positive operators in \,$\mathcal{G}\,\mathcal{B}\,(\,H\,)$\, and \,$T,\, U$\, are invertible operators in \,$\mathcal{G}\,\mathcal{B}\,(\,H\,)$.

\section{Preliminaries}
\smallskip\hspace{.6 cm}In this section, we recall some necessary definitions and theorems.

\begin{theorem}\cite{O}\label{tttt1}
Let \,$H_{\,1},\, H_{\,2}$\; be two Hilbert spaces and \;$U \,:\, H_{\,1} \,\to\, H_{\,2}$\; be a bounded linear operator with closed range \;$\mathcal{R}_{\,U}$.\;Then there exists a bounded linear operator \,$U^{\dagger} \,:\, H_{\,2} \,\to\, H_{\,1}$\, such that \,$U\,U^{\dagger}\,x \,=\, x\; \;\forall\; x \,\in\, \mathcal{R}_{\,U}$.
\end{theorem}

The operator \,$U^{\dagger}$\, defined in Theorem \ref{tttt1}, is called the pseudo-inverse of \,$U$.

\begin{theorem}(\,Douglas' factorization theorem\,)\,{\cite{Douglas}}\label{th1}
Let \;$S,\, V \,\in\, \mathcal{B}\,(\,H\,)$.\,Then the following conditions are equivalent:
\begin{description}
\item[$(i)$]$\mathcal{R}\,(\,S\,) \,\subseteq\, \mathcal{R}\,(\,V\,)$.
\item[$(ii)$]\;\;$S\, S^{\,\ast} \,\leq\, \lambda^{\,2}\; V\,V^{\,\ast}$\; for some \,$\lambda \,>\, 0$.
\item[$(iii)$]$S \,=\, V\,W$\, for some bounded linear operator \,$W$\, on \,$H$.
\end{description}
\end{theorem}

\begin{theorem}\cite{O}\label{th1.001}
The set \,$\mathcal{S}\,(\,H\,)$\; of all self-adjoint operators on \,$H$\; is a partially ordered set with respect to the partial order \,$\leq$\, which is defined as for \,$R,\,S \,\in\, \mathcal{S}\,(\,H\,)$ 
\[R \,\leq\, S \,\Leftrightarrow\, \left<\,R\,f,\, f\,\right> \,\leq\, \left<\,S\,f,\, f\,\right>\; \;\forall\; f \,\in\, H.\] 
\end{theorem}

\begin{definition}\cite{Kreyzig}
A self-adjoint operator \,$U \,:\, H_{1} \,\to\, H_{1}$\, is called positive if \,$\left<\,U\,x \,,\,  x\,\right> \,\geq\, 0$\, for all \,$x \,\in\, H_{1}$.\;In notation, we can write \,$U \,\geq\, 0$.\;A self-adjoint operator \,$V \,:\, H_{1} \,\to\, H_{1}$\, is called a square root of \,$U$\, if \,$V^{\,2} \,=\, U$.\;If, in addition \,$V \,\geq\, 0$, then \,$V$\, is called positive square root of \,$U$\, and is denoted by \,$V \,=\, U^{1 \,/\, 2}$. 
\end{definition}

\begin{theorem}\cite{Kreyzig}\label{th1.05}
The positive square root \,$V \,:\, H_{1} \,\to\, H_{1}$\, of an arbitrary positive self-adjoint operator \,$U \,:\, H_{1} \,\to\, H_{1}$\, exists and is unique.\;Further, the operator \,$V$\, commutes with every bounded linear operator on \,$H_{1}$\, which commutes with \,$U$.
\end{theorem}

In a complex Hilbert space, every bounded positive operator is self-adjoint and any two bounded positive operators can be commute with each other.

\begin{theorem}\cite{Gavruta}\label{th1.01}
Let \,$M \,\subset\, H$\; be a closed subspace and \,$T \,\in\, \mathcal{B}\,(\,H\,)$.\;Then \,$P_{\,M}\, T^{\,\ast} \,=\, P_{\,M}\,T^{\,\ast}\, P_{\,\overline{T\,M}}$.\;If \,$T$\; is an unitary operator (\,i\,.\,e \,$T^{\,\ast}\, T \,=\, I_{H}$\,), then \,$P_{\,\overline{T\,M}}\;T \,=\, T\,P_{\,M}$.
\end{theorem}

\begin{definition}\cite{Ahmadi}
Let \,$\left\{\,W_{j}\,\right\}_{ j \,\in\, J}$\, be a collection of closed subspaces of \,$H$\; and \,$\left\{\,v_{j}\,\right\}_{ j \,\in\, J}$\; be a collection of positive weights, \,$\left\{\,H_{j}\,\right\}_{ j \,\in\, J}$\, be a sequence of Hilbert spaces and let \,$\Lambda_{j} \,\in\, \mathcal{B}\,(\,H,\, H_{j}\,)$\; for each \,$j \,\in\, J$.\;Then \,$\Lambda \,=\, \{\,\left(\,W_{j},\, \Lambda_{j},\, v_{j}\,\right)\,\}_{j \,\in\, J}$\; is called a generalized fusion frame or a g-fusion frame for \,$H$\; respect to \,$\left\{\,H_{j}\,\right\}_{j \,\in\, J}$\; if there exist constants \,$0 \,<\, A \,\leq\, B \,<\, \infty$\, such that
\begin{equation}\label{eq1}
A \;\left \|\,f \,\right \|^{\,2} \,\leq\, \sum\limits_{\,j \,\in\, J}\,v_{j}^{\,2}\, \left\|\,\Lambda_{j}\,P_{\,W_{j}}\,(\,f\,) \,\right\|^{\,2} \,\leq\, B \; \left\|\, f \, \right\|^{\,2}\; \;\forall\; f \,\in\, H.
\end{equation}
The constants \,$A$\; and \,$B$\; are called the lower and upper bounds of g-fusion frame, respectively.\,If \,$A \,=\, B$\; then \,$\Lambda$\; is called tight g-fusion frame and if \;$A \,=\, B \,=\, 1$\, then we say \,$\Lambda$\; is a Parseval g-fusion frame.\;If  \,$\Lambda$\; satisfies only the right inequality of (\ref{eq1}) it is called a g-fusion Bessel sequence in \,$H$\, with bound \,$B$. 
\end{definition}

Define the space
\[l^{\,2}\left(\,\left\{\,H_{j}\,\right\}_{ j \,\in\, J}\,\right) \,=\, \left \{\,\{\,f_{\,j}\,\}_{j \,\in\, J} \,:\, f_{\,j} \;\in\; H_{j},\; \sum\limits_{\,j \,\in\, J}\, \left \|\,f_{\,j}\,\right \|^{\,2} \,<\, \infty \,\right\}\]
with inner product is given by 
\[\left<\,\{\,f_{\,j}\,\}_{ j \,\in\, J} \,,\, \{\,g_{\,j}\,\}_{ j \,\in\, J}\,\right> \;=\; \sum\limits_{\,j \,\in\, J}\, \left<\,f_{\,j} \,,\, g_{\,j}\,\right>_{H_{j}}.\]\,Clearly \,$l^{\,2}\left(\,\left\{\,H_{j}\,\right\}_{ j \,\in\, J}\,\right)$\; is a Hilbert space with the pointwise operations \cite{Sadri}. 

\begin{definition}\cite{HS}
Let \,$\left\{\,W_{j}\,\right\}_{ j \,\in\, J}$\, be a collection of closed subspaces of \,$H$\, and \,$\left\{\,v_{j}\,\right\}_{ j \,\in\, J}$\, be a collection of positive weights.\,Let \,$\left\{\,H_{j}\,\right\}_{ j \,\in\, J}$\, be a sequence of Hilbert spaces, \,$T,\, U \,\in\, \mathcal{G}\,\mathcal{B}\,(\,H\,)$\, and \,$\Lambda_{j} \,\in\, \mathcal{B}\,(\,H,\, H_{j}\,)$\, for each \,$j \,\in\, J$.\,Then the family \,$\Lambda_{T\,U} \,=\, \left\{\,\left(\,W_{j},\, \Lambda_{j},\, v_{j}\,\right)\,\right\}_{j \,\in\, J}$\, is a \,$(\,T,\,U\,)$-controlled $g$-fusion frame for \,$H$\, if there exist constants \,$0 \,<\, A \,\leq\, B \,<\, \infty$\, such that 
\begin{equation}\label{eqn1.1}
A\,\|\,f\,\|^{\,2} \,\leq\, \sum\limits_{\,j \,\in\, J}\, v^{\,2}_{j}\,\left<\,\Lambda_{j}\,P_{\,W_{j}}\,U\,f,\,  \Lambda_{j}\,P_{\,W_{j}}\,T\,f\,\right> \,\leq\, \,B\,\|\,f \,\|^{\,2}\; \;\forall\; f \,\in\, H.
\end{equation}
If \,$A \,=\, B$\, then \,$\Lambda_{T\,U}$\, is called \,$(\,T,\,U\,)$-controlled tight g-fusion frame and if \,$A \,=\, B \,=\, 1$\, then we say \,$\Lambda_{T\,U}$\, is a \,$(\,T,\,U\,)$-controlled Parseval g-fusion frame.\,If \,$\Lambda_{T\,U}$\, satisfies only the right inequality of (\ref{eqn1.1}) it is called a \,$(\,T,\,U\,)$-controlled g-fusion Bessel sequence in \,$H$.    
\end{definition}

\begin{definition}\cite{HS}
Let \,$\Lambda_{T\,U}$\, be a \,$(\,T,\,U\,)$-controlled g-fusion Bessel sequence in \,$H$\, with a bound \,$B$.\,The synthesis operator \,$T_{C} \,:\, \mathcal{K}_{\,\Lambda_{j}} \,\to\, H$\, is defined as 
\[T_{C}\,\left(\,\left\{\,v_{\,j}\,\left(\,T^{\,\ast}\,P_{\,W_{j}}\, \Lambda_{j}^{\,\ast}\,\Lambda_{j}\,P_{\,W_{j}}\,U\,\right)^{1 \,/\, 2}\,f\,\right\}_{j \,\in\, J}\,\right) \,=\, \sum\limits_{\,j \,\in\, J}\,v^{\,2}_{j}\,T^{\,\ast}\,P_{\,W_{j}}\, \Lambda_{j}^{\,\ast}\,\Lambda_{j}\,P_{\,W_{j}}\,U\,f,\]for all \,$f \,\in\, H$\, and the analysis operator \,$T^{\,\ast}_{C} \,:\, H \,\to\, \mathcal{K}_{\,\Lambda_{j}}$\,is given by 
\[T_{C}^{\,\ast}\,f \,=\,  \left\{\,v_{\,j}\,\left(\,T^{\,\ast}\,P_{\,W_{j}}\, \Lambda_{j}^{\,\ast}\,\Lambda_{j}\,P_{\,W_{j}}\,U\,\right)^{1 \,/\, 2}\,f\,\right\}_{j \,\in\, J}\; \;\forall\; f \,\in\, H,\]
where 
\[\mathcal{K}_{\,\Lambda_{j}} \,=\, \left\{\,\left\{\,v_{\,j}\,\left(\,T^{\,\ast}\,P_{\,W_{j}}\, \Lambda_{j}^{\,\ast}\,\Lambda_{j}\,P_{\,W_{j}}\,U\,\right)^{1 \,/\, 2}\,f\,\right\}_{j \,\in\, J} \,:\, f \,\in\, H\,\right\} \,\subset\, l^{\,2}\left(\,\left\{\,H_{j}\,\right\}_{ j \,\in\, J}\,\right).\]
The frame operator \,$S_{C} \,:\, H \,\to\, H$\; is defined as follows:
\[S_{C}\,f \,=\, T_{C}\,T_{C}^{\,\ast}\,f \,=\, \sum\limits_{\,j \,\in\, J}\, v_{j}^{\,2}\,T^{\,\ast}\,P_{\,W_{j}}\, \Lambda_{j}^{\,\ast}\,\Lambda_{j}\,P_{\,W_{j}}\,U\,f\; \;\forall\; f \,\in\, H\]and it is easy to verify that 
\[\left<\,S_{C}\,f,\, f\,\right> \,=\, \sum\limits_{\,j \,\in\, J}\, v^{\,2}_{j}\,\left<\,\Lambda_{j}\,P_{\,W_{j}}\,U\,f,\,  \Lambda_{j}\,P_{\,W_{j}}\,T\,f\,\right>\; \;\forall\; f \,\in\, H.\]
Furthermore, if \,$\Lambda_{T\,U}$\, is a \,$(\,T,\,U\,)$-controlled g-fusion frame with bounds \,$A$\, and \,$B$\, then \,$A\,I_{\,H} \,\leq\,S_{C} \,\leq\, B\,I_{H}$.\,Hence, \,$S_{C}$\, is bounded, invertible, self-adjoint and positive linear operator.\,It is easy to verify that \,$B^{\,-1}\,I_{H} \,\leq\, S_{C}^{\,-1} \,\leq\, A^{\,-1}\,I_{H}$.
\end{definition}

\begin{definition}\cite{GR}
Let \,$K \,\in\, \mathcal{B}\,(\,H\,)$\, and \,$\left\{\,W_{j}\,\right\}_{ j \,\in\, J}$\, be a collection of closed subspaces of \,$H$\, and \,$\left\{\,v_{j}\,\right\}_{ j \,\in\, J}$\, be a collection of positive weights.\,Let \,$\left\{\,H_{j}\,\right\}_{ j \,\in\, J}$\, be a sequence of Hilbert spaces, \,$T,\, U \,\in\, \mathcal{G}\,\mathcal{B}\,(\,H\,)$\, and \,$\Lambda_{j} \,\in\, \mathcal{B}\,(\,H,\, H_{j}\,)$\, for each \,$j \,\in\, J$.\,Then the family \,$\Lambda_{T\,U} \,=\, \left\{\,\left(\,W_{j},\, \Lambda_{j},\, v_{j}\,\right)\,\right\}_{j \,\in\, J}$\, is a \,$(\,T,\,U\,)$-controlled $K$-$g$-fusion frame for \,$H$\, if there exist constants \,$0 \,<\, A \,\leq\, B \,<\, \infty$\, such that 
\[A\,\|\,K^{\,\ast}\,f\,\|^{\,2} \,\leq\, \sum\limits_{\,j \,\in\, J}\, v^{\,2}_{j}\,\left<\,\Lambda_{j}\,P_{\,W_{j}}\,U\,f,\,  \Lambda_{j}\,P_{\,W_{j}}\,T\,f\,\right> \,\leq\, \,B\,\|\,f \,\|^{\,2}\; \;\forall\; f \,\in\, H.\]
\end{definition}

\begin{definition}\cite{MF}
Let \,$F \,:\, X \,\to\, \mathbb{H}$\; be such that for each \,$h \,\in\, H$, the mapping \,$x \,\to\, P_{\,F\,(\,x\,)}\,(\,h\,)$\; is measurable (\,i.\,e. is weakly measurable\,) and \,$v \,:\, X \,\to\, \mathbb{R}^{\,+}$\, be a measurable function and let \,$\left\{\,K_{x}\,\right\}_{x \,\in\, X}$\, be a collection of Hilbert spaces.\;For each \,$x \,\in\, X$, suppose that \,$\,\Lambda_{x} \,\in\, \mathcal{B}\,(\,F\,(\,x\,) \,,\, K_{x}\,)$.\;Then \,$\Lambda_{F} \,=\, \left\{\,\left(\,F\,(\,x\,),\, \Lambda_{x},\, v\,(\,x\,)\,\right)\,\right\}_{x \,\in\, X}$\; is called a generalized continuous fusion frame or a gc-fusion frame for \,$H$\, with respect to \,$(\,X,\, \mu\,)$\, and \,$v$, if there exists \,$0 \,<\, A \,\leq\, B \,<\, \infty$\; such that
\[A\, \|\,h\,\|^{\,2} \,\leq\, \int\limits_{\,X}\, v^{\,2}\,(\,x\,)\, \left\|\,\Lambda_{x}\,P_{\,F\,(\,x\,)}\,(\,h\,)\,\right\|^{\,2}\,d\mu \,\leq\, B\, \|\,h\,\|^{\,2}\;\; \;\forall\, h \,\in\, H,\]where \,$P_{\,F\,(\,x\,)}$\, is the orthogonal projection onto the subspace \,$F\,(\,x\,)$.\;$\Lambda_{F}$\, is called a tight gc-fusion frame for \,$H$\, if \,$A \,=\, B$\, and Parseval if \,$A \,=\, B \,=\, 1$.\;If we have only the upper bound, we call \,$\Lambda_{F}$\, is a Bessel gc-fusion mapping for \,$H$.
\end{definition}

Let \,$K \,=\, \oplus_{x \,\in\, X}\,K_{x}$\, and \,$L^{\,2}\left(\,X,\, K\,\right)$\, be a collection of all measurable functions \,$\varphi \,:\, X \,\to\, K$\, such that for each \,$x \,\in\, X,\, \varphi\,(\,x\,) \,\in\, K_{x}$\; and \,$\int\limits_{\,X}\,\left\|\,\varphi\,(\,x\,)\,\right\|^{\,2}\,d\mu \,<\, \infty$. It can be verified that \,$L^{\,2}\left(\,X,\, K\,\right)$\; is a Hilbert space with inner product given by
\[\left<\,\phi,\, \varphi\,\right> \,=\, \int\limits_{\,X}\, \left<\,\phi\,(\,x\,),\, \varphi\,(\,x\,)\,\right>\,d\mu\] for \,$\phi,\, \varphi \,\in\, L^{\,2}\left(\,X,\, K\,\right)$.    

\begin{definition}\cite{MF}
Let \,$\Lambda_{F} \,=\, \left\{\,\left(\,F\,(\,x\,),\, \Lambda_{x},\, v\,(\,x\,)\,\right)\,\right\}_{x \,\in\, X}$\, be a Bessel gc-fusion mapping for \,$H$.\;Then the gc-fusion pre-frame operator or synthesis operator \,$T_{g\,F} \,:\, L^{\,2}\left(\,X,\, K\,\right) \,\to\, H$\, is defined by
\[\left<\,T_{g\,F}\,(\,\varphi\,),\, h\,\right> \,=\, \int\limits_{\,X}\, v\,(\,x\,)\,\left<\,P_{\,F\,(\,x\,)}\,\Lambda^{\,\ast}_{x}\,\left(\,\varphi\,(\,x\,)\,\right),\, h\,\right>\]where \,$\varphi \,\in\, L^{\,2}\left(\,X,\, K\,\right)$\, and \,$h \,\in\, H$.\,$T_{g\,F}$\, is a bounded linear mapping and its adjoint operator is given by  
\[T^{\,\ast}_{g\,F} \,:\, H \,\to\, L^{\,2}\left(\,X,\, K\,\right),\; T^{\,\ast}_{g\,F}\,(\,h\,) \,=\, \left\{\,v\,(\,x\,)\, \Lambda_{x}\,P_{\,F\,(\,x\,)}\,(\,h\,)\,\right\}_{x \,\in\, X},\; h \,\in\, H\]and \,$S_{g\,F} \,=\, T_{g\,F}\,T^{\,\ast}_{g\,F}$\, is called gc-fusion frame operator.\;Thus, for each \,$f,\, h \,\in\, H$,
\[\left<\,S_{g\,F}\,(\,f\,),\, h\,\right> \,=\, \int\limits_{\,X}\, v^{\,2}\,(\,x\,)\,\left<\,P_{F\,(\,x\,)}\,\Lambda^{\,\ast}_{x}\,\Lambda_{x}\,P_{F\,(\,x\,)}\,f,\; h\,\right>.\]The operator \,$S_{g\,F}$\, is bounded, self-adjoint, positive and invertible operator on \,$H$\,.
\end{definition}

\section{Continuous controlled $K$-$g$-fusion frame }

\smallskip\hspace{.6 cm}In this section, continuous version of controlled $K$-$g$-fusion frame for \,$H$\, is presented.\,We expand some of the recent results of controlled \,$K$-$g$-fusion frame to continuous controlled $K$-$g$-fusion frame.

\begin{definition}
Let \,$K \,\in\, \mathcal{B}\,(\,H\,)$\, and \,$F \,:\, X \,\to\, \mathbb{H}$\, be a mapping, \,$v \,:\, X \,\to\, \mathbb{R}^{\,+}$\, be a measurable function and \,$\left\{\,K_{x}\,\right\}_{x \,\in\, X}$\, be a collection of Hilbert spaces.\;For each \,$x \,\in\, X$, suppose that \,$\,\Lambda_{x} \,\in\, \mathcal{B}\,(\,F\,(\,x\,),\, K_{x}\,)$\, and \,$T,\, U \,\in\, \mathcal{G}\,\mathcal{B}^{\,+}\,(\,H\,)$.\,Then \,$\Lambda_{T\,U} \,=\, \left\{\,\left(\,F\,(\,x\,),\, \Lambda_{x},\, v\,(\,x\,)\,\right)\,\right\}_{x \,\in\, X}$\, is called a continuous \,$(\,T,\,U\,)$-controlled $K$-$g$-fusion frame for \,$H$\, with respect to \,$(\,X,\, \mu\,)$\, and \,$v$, if
\begin{description}
\item[$(i)$]for each \,$f \,\in\, H$, the mapping \,$x \,\to\, P_{F\,(\,x\,)}\,(\,f\,)$\; is measurable (\,i.\,e. is weakly measurable\,).
\item[$(ii)$]there exist constants \,$0 \,<\, A \,\leq\, B \,<\, \infty$\, such that
\end{description}
\begin{equation}\label{eq1}
A\,\left\|\,K^{\,\ast}\,f\,\right\|^{\,2} \,\leq\, \int\limits_{\,X}\,v^{\,2}\,(\,x\,)\,\left<\,\Lambda_{x}\,P_{\,F\,(\,x\,)}\,U\,f,\, \Lambda_{x}\,P_{\,F\,(\,x\,)}\,T\,f\,\right>\,d\mu_{x} \,\leq\, B\,\|\,f\,\|^{\,2},
\end{equation}
for all \,$f \,\in\, H$, where \,$P_{\,F\,(\,x\,)}$\, is the orthogonal projection onto the subspace \,$F\,(\,x\,)$.\,The constants \,$A,\,B$\, are called the frame bounds.
\end{definition}
Furthermore,
\begin{description}
\item[$(I)$]if only the right inequality of (\ref{eq1}) holds then \,$\Lambda_{T\,U}$\, is called a continuous \,$(\,T,\,U\,)$-controlled $K$-$g$-fusion Bessel family for \,$H$.
\item[$(II)$]if \,$T \,=\, I_{H}$\, then \,$\Lambda_{T\,U}$\, is called a continuous \,$(\,I_{H},\, U\,)$-controlled $K$-$g$-fusion frame for \,$H$. 
\item[$(III)$]if \,$T \,=\, U \,=\, I_{H}$\, then \,$\Lambda_{T\,U}$\, is called a continuous $K$-$g$-fusion frame for \,$H$.
\item[$(IV)$]if \,$K \,=\, I_{H}$\, then \,$\Lambda_{T\,U}$\, is called a continuous \,$(\,T,\,U\,)$-controlled $g$-fusion frame for \,$H$.  
\end{description}
 
\begin{remark}
If the measure space \,$X \,=\, \mathbb{N}$\, and \,$\mu$\, is the counting measure then a continuous \,$(\,T,\,U\,)$-controlled $g$-fusion frame will be the discrete \,$(\,T,\,U\,)$-controlled $g$-fusion frame.   
\end{remark} 

\begin{proposition}\label{pp1}
Let \,$\Lambda_{T\,U}$\, be a continuous \,$(\,T,\,U\,)$-controlled $g$-fusion Bessel family for \,$H$\, with bound \,$B$.\,Then there exists a unique bounded linear operator \,$S_{C} \,:\, H \,\to\,\ H$\, such that
\[\left<\,S_{C}\,f,\, g\,\right> \,=\, \int\limits_{\,X}\,v^{\,2}\,(\,x\,)\,\left<\,T^{\,\ast}\,P_{F\,(\,x\,)}\,\Lambda_{x}^{\,\ast}\,\Lambda_{x}\,P_{F\,(\,x\,)}\,U\,f,\, g\,\right>\,d\mu_{x}\; \;\forall\; f,\, g \,\in\, H.\]
Furthermore, if \,$\Lambda_{T\,U}$\, is a continuous \,$(\,T,\,U\,)$-controlled $K$-$g$-fusion frame for \,$H$\, then \,$A\,K\,K^{\,\ast} \,\leq\, S_{C} \,\leq\, B\,I_{H}$. 
\end{proposition}

\begin{proof}
Proof of this proposition is directly follows from the proposition 3.3 of \cite{GP}.\\
  
Furthermore, if \,$\Lambda_{T\,U}$\, is a continuous \,$(\,T,\,U\,)$-controlled $K$-$g$-fusion frame for \,$H$\, then by (\ref{eq1}) it  is easy to verify that \,$A\,K\,K^{\,\ast} \,\leq\, S_{C} \,\leq\, B\,I_{H}$.  
\end{proof}
 
The operator defined in the Proposition \ref{pp1} is called the frame operator for \,$\Lambda_{T\,U}$.   
 
\begin{definition}
Let \,$\Lambda_{T\,U}$\, be a continuous \,$(\,T,\,U\,)$-controlled $g$-fusion Bessel family for \,$H$.\,Then the bounded linear operator \,$T_{C} \,:\, L^{\,2}\left(\,X,\, K\,\right) \,\to\, H$\, defined by 
\[\left<\,T_{C}\,\Phi,\, g\,\right> \,=\, \int\limits_{\,X}\,v^{\,2}\,(\,x\,)\,\left<\,T^{\,\ast}\,P_{F\,(\,x\,)}\,\Lambda_{x}^{\,\ast}\,\Lambda_{x}\,P_{F\,(\,x\,)}\,U\,f,\, g\,\right>\,d\mu_{x},\]
where for all \,$f \,\in\, H$, \,$\Phi \,=\, \left\{\,v\,(\,x\,)\,\left(\,T^{\,\ast}\,P_{F\,(\,x\,)}\,\Lambda_{x}^{\,\ast}\,\Lambda_{x}\,P_{F\,(\,x\,)}\,U\,\right)^{1 \,/\, 2}\,f\,\right\}_{x \,\in\, X}$\, and \,$g \,\in\, H$\, is called synthesis operator and its adjoint operator described by 
\[T_{C}^{\,\ast}\,g \,=\, \left\{\,v\,(\,x\,)\,\left(\,T^{\,\ast}\,P_{F\,(\,x\,)}\,\Lambda_{x}^{\,\ast}\,\Lambda_{x}\,P_{F\,(\,x\,)}\,U\,\right)^{1 \,/\, 2}\,g\,\right\}_{x \,\in\, X}.\]
is called analysis operator.
\end{definition}
 
Next we will see that continuous controlled \,$g$-fusion Bessel families for \,$H$\, becomes continuous controlled \,$g$-fusion frames for \,$H$\, under some sufficient conditions. 

\begin{theorem}
Let the families \,$\Lambda_{T\,U} \,=\, \left\{\,\left(\,F\,(\,x\,),\, \Lambda_{x},\, v\,(\,x\,)\,\right)\,\right\}_{x \,\in\, X}$\, and \,$\Gamma_{T\,U} \,=\, \left\{\,\left(\,F\,(\,x\,),\, \Gamma_{x},\, v\,(\,x\,)\,\right)\,\right\}_{x \,\in\, X}$\, be two continuous \,$(\,T,\,U\,)$-controlled $g$-fusion Bessel families for \,$H$\, with bounds \,$B$\, and \,$D$, respectively.\,Suppose that \,$T_{C}$\, and \,$T_{C}^{\,\prime}$\, be their synthesis operators such that \,$T_{C}^{\,\prime}\,T^{\,\ast}_{C} \,=\, K^{\,\ast}$.\,Then \,$\Lambda_{T\,U}$\, and \,$\Gamma_{T\,U}$\, are continuous \,$(\,T,\,U\,)$-controlled $K$-$g$-fusion frame and continuous \,$(\,T,\,U\,)$-controlled $K^{\,\ast}$-$g$-fusion frame for \,$H$, respectively.  
\end{theorem}

\begin{proof}
For each \,$f \,\in\, H$, we have
\begin{align*}
\left\|\,K^{\,\ast}\,f\,\right\|^{\,4} &\,= \left<\,K^{\,\ast}\,f,\, K^{\,\ast}\,f\,\right>^{\,2} \,=\, \left<\,T^{\,\ast}_{C}\,f,\, \left(\,T_{C}^{\,\prime}\,\right)^{\,\ast}K^{\,\ast}\,f\,\right>^{\,2} \,\leq\, \left\|\,T^{\,\ast}_{C}\,f\,\right\|^{\,2}\,\left\|\,\left(\,T_{C}^{\,\prime}\,\right)^{\,\ast}K^{\,\ast}\,f\,\right\|^{\,2}\\
&=\,\int\limits_{\,X}\,v^{\,2}\,(\,x\,)\,\left<\,\Lambda_{x}\,P_{F\,(\,x\,)}\,U\,f,\, \Lambda_{x}\,P_{F\,(\,x\,)}\,T\,f\,\right>\,d\mu_{x}\;\times\\
&\hspace{1cm}\int\limits_{\,X}\,v^{\,2}\,(\,x\,)\,\left<\,\Gamma_{x}\,P_{F\,(\,x\,)}\,U\,K^{\,\ast}\,f,\, \Gamma_{x}\,P_{F\,(\,x\,)}\,T\,K^{\,\ast}\,f\,\right>\,d\mu_{x}\,\\
&\leq\, D\,\left\|\,K^{\,\ast}\,f\,\right\|^{\,2}\,\int\limits_{\,X}\,v^{\,2}\,(\,x\,)\,\left<\,\Lambda_{x}\,P_{F\,(\,x\,)}\,U\,f,\, \Lambda_{x}\,P_{F\,(\,x\,)}\,T\,f\,\right>\,d\mu_{x}\\
&\Rightarrow\, \dfrac{1}{D}\,\left\|\,K^{\,\ast}\,f\,\right\|^{\,2} \,\leq\, \int\limits_{\,X}\,v^{\,2}\,(\,x\,)\,\left<\,\Lambda_{x}\,P_{F\,(\,x\,)}\,U\,f,\, \Lambda_{x}\,P_{F\,(\,x\,)}\,T\,f\,\right>\,d\mu_{x}.  
\end{align*}
This shows that \,$\Lambda_{T\,U}$\, is a continuous \,$(\,T,\,U\,)$-controlled $K$-$g$-fusion frame for \,$H$\, with bounds \,$1 \,/\, D$\, and \,$B$.\,Similarly, it can be shown that \,$\Gamma_{T\,U}$\, is a continuous \,$(\,T,\,U\,)$-controlled $K^{\,\ast}$-$g$-fusion frame for \,$H$.  
\end{proof}

In the following theorem, we will see that any continuous controlled $K$-$g$-fusion frame is a continuous $K$-$g$-fusion frame and conversely any continuous $K$-$g$-fusion frame is a continuous controlled $K$-$g$-fusion frame under some sufficient conditions. 

\begin{theorem}\label{th1.3}
Let \,$T,\, U \,\in\, \mathcal{G}\,\mathcal{B}^{\,+}\,(\,H\,)$\, and \,$S_{g\,F}\,T \,=\, T\,S_{g\,F}$.\,If the operator \,$K$\, commutes with \,$T$\, and \,$U$\, then \,$\Lambda_{T\,U}$\, is a continuous \,$(\,T,\,U\,)$-controlled $K$-$g$-fusion frame for \,$H$\, if and only if \,$\Lambda_{T\,U}$\, is a continuous $K$-$g$-fusion frame for \,$H$, where \,$S_{g\,F}$\, is the continuous \,$g$-fusion frame operator defined by
\[\left<\,S_{g\,F}\,f,\, f\,\right> \,=\, \int\limits_{\,X}\,v^{\,2}\,(\,x\,)\,\left<\,P_{F\,(\,x\,)}\,\Lambda^{\,\ast}_{x}\,\Lambda_{x}\,P_{F\,(\,x\,)}\,f,\, f\,\right>\,d\mu_{x},\; \;f \,\in\, H.\] 
\end{theorem}

\begin{proof}
First we suppose that \,$\Lambda_{T\,U}$\, is a continuous $K$-$g$-fusion frame for \,$H$\, with bounds \,$A$\, and \,$B$.\,Then for each \,$f \,\in\, H$, we have
\[A\,\left\|\,K^{\,\ast}\,f\,\right\|^{\,2} \,\leq\, \int\limits_{\,X}\,v^{\,2}\,(\,x\,)\,\left\|\,\Lambda_{x}\,P_{F\,(\,x\,)}\,f\,\right\|^{\,2}\,d\mu_{x} \,\leq\, B\,\|\,f\,\|^{\,2}.\]
Now according to the Lemma 3.10 of \cite{NA}, we can deduced that 
\[m\,m^{\,\prime}\,A\,K\,K^{\,\ast} \,\leq\, T\,S_{g\,F}\,U \,\leq\, M\,M^{\prime}\,B\,I_{H},\] where \,$m,\,m^{\,\prime}$\, and \,$M,\,M^{\prime}$\, are positive constants.\,Then for each \,$f \,\in\, H$, we have
\begin{align*}
&m\,m^{\,\prime}\,A\,\left\|\,K^{\,\ast}\,f\,\right\|^{\,2} \,\leq\, \int\limits_{\,X}\,v^{\,2}\,(\,x\,)\,\left<\,T\,P_{F\,(\,x\,)}\,\Lambda^{\,\ast}_{x}\,\Lambda_{x}\,P_{F\,(\,x\,)}\,U\,f,\, f\,\right>\,d\mu_{x} \,\leq\, M\,M^{\,\prime}\,B\,\|\,f \,\|^{\,2}\\
&\Rightarrow\, m\,m^{\,\prime}\,A\,\left\|\,K^{\,\ast}\,f\,\right\|^{\,2} \,\leq\, \int\limits_{\,X}\,v^{\,2}\,(\,x\,)\,\left<\,\Lambda_{x}\,P_{F\,(\,x\,)}\,U\,f,\, \Lambda_{x}\,P_{F\,(\,x\,)}\,T\,f\,\right>\,d\mu_{x} \,\leq\, M\,M^{\,\prime}\,B\,\|\,f \,\|^{\,2}.
\end{align*}
Hence, \,$\Lambda_{T\,U}$\, is a continuous \,$(\,T,\,U\,)$-controlled $K$-$g$-fusion frame for \,$H$.\\ 

Conversely, suppose that \,$\Lambda_{T\,U}$\, is a continuous \,$(\,T,\,U\,)$-controlled $K$-$g$-fusion frame for \,$H$\, with bounds \,$A$\, and \,$B$.\,Now, for each \,$f \,\in\, H$, we have
\begin{align*}
&A\,\left\|\,K^{\,\ast}\,f\,\right\|^{\,2} \,=\, A\,\left\|\,(\,T\,U\,)^{1 \,/\, 2}\,(\,T\,U\,)^{\,-\, 1 \,/\, 2}\,K^{\,\ast}\,f\,\right\|^{\,2} \,=\, A\,\left\|\,(\,T\,U\,)^{1 \,/\, 2}\,K^{\,\ast}\,(\,T\,U\,)^{\,-\, 1 \,/\, 2}\,f\,\right\|^{\,2}\\ 
&\leq\, \left\|\,(\,T\,U\,)^{1 \,/\, 2}\,\right\|^{\,2}\, \int\limits_{\,X}\,v^{\,2}\,(\,x\,)\,\left<\,\Lambda_{x}\,P_{F\,(\,x\,)}\,U\,(\,T\,U\,)^{\,-\, 1 \,/\, 2}\,f,\, \Lambda_{x}\,P_{F\,(\,x\,)}\,T\,(\,T\,U\,)^{\,-\, 1 \,/\, 2}\,f\,\right>\,d\mu_{x} \\
&=\, \left\|\,(\,T\,U\,)^{1 \,/\, 2}\,\right\|^{\,2}\,\int\limits_{\,X}\,v^{\,2}\,(\,x\,)\,\left<\,\Lambda_{x}\,P_{F\,(\,x\,)}\,U^{1 \,/\, 2}\,T^{\,-\, 1 \,/\, 2}\,f,\, \Lambda_{x}\,P_{F\,(\,x\,)}\,T^{1 \,/\, 2}\,U^{\,-\, 1 \,/\, 2}\,f\,\right>\,d\mu_{x}\\
&=\,\left\|\,(\,T\,U\,)^{1 \,/\, 2}\,\right\|^{\,2}\,\int\limits_{\,X}\,v^{\,2}\,(\,x\,)\,\left<\,U^{\,-\, 1 \,/\, 2}\,T^{1 \,/\, 2}\,P_{F\,(\,x\,)}\,\Lambda^{\,\ast}_{x}\,\Lambda_{x}\,P_{F\,(\,x\,)}\,U^{1 \,/\, 2}\,T^{\,-\, 1 \,/\, 2}\,f,\, f\,\right>\,d\mu_{x}\\ 
&= \left\|\,(\,T\,U\,)^{1 \,/\, 2}\,\right\|^{\,2} \left<\,U^{\,-\, 1 \,/\, 2}\,T^{1 \,/\, 2}\,S_{g\,F}\,U^{1 \,/\, 2}\,T^{\,-\, 1 \,/\, 2}\,f,\, f\,\right> = \left\|\,(\,T\,U\,)^{1 \,/\, 2}\,\right\|^{\,2} \left<\,S_{g\,F}\,f,\, f\,\right>\\  
&=\, \left\|\,(\,T\,U\,)^{1 \,/\, 2}\,\right\|^{\,2}\,\int\limits_{\,X}\,v^{\,2}\,(\,x\,)\,\left<\,P_{F\,(\,x\,)}\,\Lambda^{\,\ast}_{x}\,\Lambda_{x}\,P_{F\,(\,x\,)}\,f,\, f\,\right>\,d\mu_{x}\\ 
&\Rightarrow\, \dfrac{A}{\left\|\,(\,T\,U\,)^{1 \,/\, 2}\,\right\|^{\,2}}\,\left\|\,K^{\,\ast}\,f\,\right\|^{\,2} \,\leq\, \int\limits_{\,X}\,v^{\,2}\,(\,x\,)\,\left\|\,\Lambda_{x}\,P_{F\,(\,x\,)}\,f\,\right\|^{\,2}\,d\mu_{x}.
\end{align*}
On the other hand, it is easy to verify that
\begin{align*}
&\int\limits_{\,X}\,v^{\,2}\,(\,x\,)\,\left\|\,\Lambda_{x}\,P_{F\,(\,x\,)}\,f\,\right\|^{\,2}\,d\mu_{x} \,=\, \left<\,(\,T\,U\,)^{\,-\, 1 \,/\, 2}\,(\,T\,U\,)^{1 \,/\, 2}\,S_{g\,F}\,f,\, f\,\right>\\
&= \left<\,(\,T\,U\,)^{1 \,/\, 2}\,S_{g\,F}\,f,\, (\,T\,U\,)^{\,-\, 1 \,/\, 2}\,f\,\right> = \left<\,S_{g\,F}\,(\,T\,U\,)\,(\,T\,U\,)^{\,-\, 1 \,/\, 2}\,f,\, (\,T\,U\,)^{\,-\, 1 \,/\, 2}\,f\,\right>\\
&=\,\left<\,T\,S_{g\,F}\,U\,(\,T\,U\,)^{\,-\, 1 \,/\, 2}\,f,\, (\,T\,U\,)^{\,-\, 1 \,/\, 2}\,f\,\right> =\left<\,S_{C}\,(\,T\,U\,)^{\,-\, 1 \,/\, 2}\,f,\, (\,T\,U\,)^{\,-\, 1 \,/\, 2}\,f\,\right>\\
&\,\leq\, B\,\left\|\,(\,T\,U\,)^{\,-\, 1 \,/\, 2}\,\right\|^{\,2}\,\|\,f\,\|^{\,2}.
\end{align*} 
Thus, \,$\Lambda_{T\,U}$\, is a continuous $K$-$g$-fusion frame for \,$H$.\,This completes the proof. 
\end{proof}

In the next two theorems, we will construct new type of continuous controlled $g$-fusion frame from a given continuous controlled \,$K$-$g$-fusion frame by using a invertible bounded linear operator.

\begin{theorem}\label{ttt1}
Let \,$\Lambda_{T\,U}$\, be a continuous \,$(\,T,\,U\,)$-controlled $K$-$g$-fusion frame for \,$H$\, with bounds \,$A,\, B$\, and \,$V \,\in\, \mathcal{B}\,(\,H\,)$\, be an invertible operator on \,$H$\, such that \,$V^{\,\ast}$\, commutes with \,$T$\, and \,$U$.\,Then \,$\Gamma_{T\,U} \,=\, \left\{\,\left(\,V\,F\,(\,x\,),\, \Lambda_{x}\,P_{F\,(\,x\,)}\,V^{\,\ast},\, v\,(\,x\,)\,\right)\,\right\}_{x \,\in\, X}$\, is a continuous \,$(\,T,\,U\,)$-controlled $V\,K\,V^{\,\ast}$-$g$-fusion frame for \,$H$.  
\end{theorem} 
 
\begin{proof}
Since \,$P_{\,F\,(\,x\,)}\,V^{\,\ast} \,=\, P_{\,F\,(\,x\,)}\,V^{\,\ast}\,P_{\,V\,F\,(\,x\,)}$\, for all \,$x \,\in\, X$, the mapping \,$f \,\to\, P_{\,V\,F\,(\,x\,)}\,f,\, \,f \,\in\, H$\, is weakly measurable.\,Now, for each \,$f \,\in\, H$, using Theorem \ref{th1.01}, we have
\begin{align}
&\int\limits_{\,X}\,v^{\,2}\,(\,x\,)\,\left<\,\Lambda_{x}\,P_{\,F\,(\,x\,)}\,V^{\,\ast}\,P_{\,V\,F\,(\,x\,)}\,U\,f,\, \Lambda_{x}\,P_{F\,(\,x\,)}\,V^{\,\ast}\,P_{V\,F\,(\,x\,)}\,T\,f\,\right>\,d\mu_{x}\nonumber\\
&=\,\int\limits_{\,X}\,v^{\,2}\,(\,x\,)\,\left<\,\Lambda_{x}\,P_{F\,(\,x\,)}\,V^{\,\ast}\,U\,f,\, \Lambda_{x}\,P_{F\,(\,x\,)}\,V^{\,\ast}\,T\,f\,\right>\,d\mu_{x}\nonumber\\ 
&=\,\int\limits_{\,X}\,v^{\,2}\,(\,x\,)\,\left<\,\Lambda_{x}\,P_{F\,(\,x\,)}\,U\,V^{\,\ast}\,f,\, \Lambda_{x}\,P_{F\,(\,x\,)}\,T\,V^{\,\ast}\,f\,\right>\,d\mu_{x}\nonumber\\
&\leq\, B\,\left\|\,V^{\,\ast}\,f\,\right\|^{\,2} \,\leq\, B\,\|\,V\,\|^{\,2}\,\|\,f\,\|^{\,2}.\nonumber
\end{align}
On the other hand, for each \,$f \,\in\, H$, we get
\begin{align*}
&\dfrac{A}{\|\,V\,\|^{\,2}}\, \left\|\,\left(\,V\,K\,V^{\,\ast}\,\right)^{\,\ast}\,f\,\right\|^{\,2} \,=\, \dfrac{A}{\|\,V\,\|^{\,2}}\, \left\|\,V\,K^{\,\ast}\,V^{\,\ast}\,f\,\right\|^{\,2} \,\leq\, A\, \left\|\,K^{\,\ast}\,V^{\,\ast}\,f\,\right\|^{\,2}\\
&=\,\int\limits_{\,X}\,v^{\,2}\,(\,x\,)\,\left<\,\Lambda_{x}\,P_{F\,(\,x\,)}\,U\,V^{\,\ast}\,f,\, \Lambda_{x}\,P_{F\,(\,x\,)}\,T\,V^{\,\ast}\,f\,\right>\,d\mu_{x}\\
&=\,\int\limits_{\,X}\,v^{\,2}\,(\,x\,)\,\left<\,\Lambda_{x}\,P_{F\,(\,x\,)}\,V^{\,\ast}\,U\,f,\, \Lambda_{x}\,P_{F\,(\,x\,)}\,V^{\,\ast}\,T\,f\,\right>\,d\mu_{x}\\ 
&=\,\int\limits_{\,X}\,v^{\,2}\,(\,x\,)\,\left<\,\Lambda_{x}\,P_{F\,(\,x\,)}\,V^{\,\ast}\,P_{V\,F\,(\,x\,)}\,U\,f,\, \Lambda_{x}\,P_{F\,(\,x\,)}\,V^{\,\ast}\,P_{V\,F\,(\,x\,)}\,T\,f\,\right>\,d\mu_{x}\\
\end{align*}
Thus, \,$\Gamma_{T\,U}$\, is a continuous \,$(\,T,\,U\,)$-controlled $V\,K\,V^{\,\ast}$-$g$-fusion frame for \,$H$\, with bounds \,$A \,/\, \|\,V\,\|^{\,2}$\, and \,$B\,\|\,V\,\|^{\,2}$.
\end{proof}

\begin{theorem}
Let \,$V \,\in\, \mathcal{B}\,(\,H\,)$\, be an invertible operator on \,$H$, \,$(\,V^{\,-\, 1}\,)^{\,\ast}$\, commutes with \,$T$\, and \,$U$.\,Suppose \,$\Gamma_{T\,U} \,=\, \left\{\,\left(\,V\,F\,(\,x\,),\, \Lambda_{x}\,P_{F\,(\,x\,)}\,V^{\,\ast},\, v\,(\,x\,)\,\right)\,\right\}_{x \,\in\, X}$\, is a continuous \,$(\,T,\,U\,)$-controlled $K$-$g$-fusion frame for \,$H$, for some \,$K \,\in\, \mathcal{B}\,(\,H\,)$.\,Then \,$\Lambda_{T\,U}$\, is a continuous \,$(\,T,\,U\,)$-controlled $V^{\,-\, 1}\,K\,V$-$g$-fusion frame for \,$H$. 
\end{theorem}

\begin{proof}
Since \,$\Gamma_{T\,U}$\, is a continuous \,$(\,T,\,U\,)$-controlled\,$K$-$g$-fusion frame for \,$H$, for each \,$f \,\in\, H$, there exist constants \,$A,\, B \,>\, 0$\, such that
\begin{align}
A \,\left \|\,K^{\,\ast}\,f \,\right \|^{\,2} &\,\leq\, \int\limits_{\,X}\,v^{\,2}\,(\,x\,)\,\left<\,\Lambda_{x}\,P_{F\,(\,x\,)}\,V^{\,\ast}\,P_{V\,F\,(\,x\,)}\,U\,f,\, \Lambda_{x}\,P_{F\,(\,x\,)}\,V^{\,\ast}\,P_{V\,F\,(\,x\,)}\,T\,f\,\right>\,d\mu_{x}\nonumber\\
&\leq\, B \,\left\|\,f\, \right\|^{\,2}.\label{eq1.3}
\end{align}
Now, for each \,$f \,\in\, H$, using Theorem \ref{th1.01}, we have
\begin{align*}
&\dfrac{A}{\|\,V\,\|^{\,2}}\,\left \|\,\left(\,V^{\,-\, 1}\,K\,V\,\right)^{\,\ast}\,f\,\right \|^{\,2} \,=\, \dfrac{A}{\|\,V\,\|^{\,2}}\,\left\|\,V^{\,\ast}\,K^{\,\ast}\,(\,V^{\,-\, 1}\,)^{\,\ast}\,f\,\right\|^{\,2}\\
& \,\leq\, A\;\left\|\,K^{\,\ast}\,\left(\,V^{\,-\, 1}\,\right)^{\,\ast}\,f\,\right\|^{\,2}\\
&\leq\, \int\limits_{\,X}\,v^{\,2}\,(\,x\,)\,\left<\,\Lambda_{x}\,P_{F\,(\,x\,)}\,V^{\,\ast}\,P_{V\,F\,(\,x\,)}\,U\,\left(\,V^{\,-\, 1}\,\right)^{\,\ast}\,f,\, \Lambda_{x}\,P_{F\,(\,x\,)}\,V^{\,\ast}\,P_{V\,F\,(\,x\,)}\,T\,\left(\,V^{\,-\, 1}\,\right)^{\,\ast}\,f\,\right>\,d\mu_{x}\\
&= \int\limits_{\,X}\,v^{\,2}\,(\,x\,)\,\left<\,\Lambda_{x}\,P_{F\,(\,x\,)}\,V^{\,\ast}\,U\,\left(\,V^{\,-\, 1}\,\right)^{\,\ast}\,f,\, \Lambda_{x}\,P_{F\,(\,x\,)}\,V^{\,\ast}\,T\,\left(\,V^{\,-\, 1}\,\right)^{\,\ast}\,f\,\right>\,d\mu_{x}\\
&=\int\limits_{\,X}\,v^{\,2}\,(\,x\,)\,\left<\,\Lambda_{x}\,P_{F\,(\,x\,)}\,V^{\,\ast}\,\left(\,V^{\,-\, 1}\,\right)^{\,\ast}\,U\,f,\, \Lambda_{x}\,P_{F\,(\,x\,)}\,V^{\,\ast}\,\left(\,V^{\,-\, 1}\,\right)^{\,\ast}\,T\,f\,\right>\,d\mu_{x}\\
&=\int\limits_{\,X}\,v^{\,2}\,(\,x\,)\,\left<\,\Lambda_{x}\,P_{F\,(\,x\,)}\,U\,f,\, \Lambda_{x}\,P_{F\,(\,x\,)}\,T\,f\,\right>\,d\mu_{x}.
\end{align*}
On the other hand, for each \,$f \,\in\, H$, we have
\begin{align*}
&\int\limits_{\,X}\,v^{\,2}\,(\,x\,)\,\left<\,\Lambda_{x}\,P_{F\,(\,x\,)}\,U\,f,\, \Lambda_{x}\,P_{F\,(\,x\,)}\,T\,f\,\right>\,d\mu_{x}\\
&= \int\limits_{\,X}\,v^{\,2}\,(\,x\,)\,\left<\,\Lambda_{x}\,P_{F\,(\,x\,)}\,V^{\,\ast}\,U\,\left(\,V^{\,-\, 1}\,\right)^{\,\ast}\,f,\, \Lambda_{x}\,P_{F\,(\,x\,)}\,V^{\,\ast}\,T\,\left(\,V^{\,-\, 1}\,\right)^{\,\ast}\,f\,\right>\,d\mu_{x}\\
&= \int\limits_{\,X}\,v^{\,2}\,(\,x\,)\,\left<\,\Lambda_{x}\,P_{F\,(\,x\,)}\,V^{\,\ast}\,P_{V\,F\,(\,x\,)}\,U\,\left(\,V^{\,-\, 1}\,\right)^{\,\ast}\,f,\, \Lambda_{x}\,P_{F\,(\,x\,)}\,V^{\,\ast}\,P_{V\,F\,(\,x\,)}\,T\,\left(\,V^{\,-\, 1}\,\right)^{\,\ast}\,f\,\right>\,d\mu_{x}\\
&\leq\, B \; \left\|\,\left(\,V^{\,-\, 1}\,\right)^{\,\ast}\,f\, \right\|^{\,2}\,\leq\, B\; \left\|\,V^{\,-\, 1}\,\right\|^{\,2}\,\|\,f\,\|^{\,2}\; [\;\text{by (\ref{eq1.3})}\;].
\end{align*}
Thus, \,$\Lambda_{T\,U}$\, is a continuous \,$(\,T,\,U\,)$-controlled $V^{\,-\, 1}\,K\,V$-$g$-fusion frame for \,$H$.
\end{proof}

In the following theorem, we will see that every continuous controlled \,$g$-fusion frame is a continuous controlled \,$K$-$g$-fusion frame and the converse is also true under some condition.

\begin{theorem}
Let \,$K \,\in\, \mathcal{B}\,(\,H\,)$.\,Then 
\begin{description}
\item[$(i)$]Every continuous \,$(\,T,\,U\,)$-controlled $g$-fusion frame is a continuous \,$(\,T,\,U\,)$-controlled $K$-$g$-fusion frame. 
\item[$(ii)$]If \,$\mathcal{R}\,(\,K\,)$\, is closed, every continuous \,$(\,T,\,U\,)$-controlled $K$-$g$-fusion frame is a continuous \,$(\,T,\,U\,)$-controlled $g$-fusion frame for \,$\mathcal{R}\,(\,K\,)$.  
\end{description}
\end{theorem}

\begin{proof}$(i)$
Let \,$\Lambda_{T\,U}$\, be a continuous \,$(\,T,\,U\,)$-controlled $g$-fusion frame for \,$H$\, with bounds \,$A$\, and \,$B$.\,Then for each \,$f \,\in\, H$, we have
\begin{align*}
&\dfrac{A}{\|\,K\,\|^{\,2}}\,\left \|\,K^{\,\ast}\,f \,\right \|^{\,2} \,\leq\, A\,\|\,f\,\|^{\,2}\\
& \,\leq\, \int\limits_{\,X}\,v^{\,2}\,(\,x\,)\,\left<\,\Lambda_{x}\,P_{F\,(\,x\,)}\,U\,f,\, \Lambda_{x}\,P_{F\,(\,x\,)}\,T\,f\,\right>\,d\mu_{x} \,\leq\, \,B\,\|\,f \,\|^{\,2}.
\end{align*} 
Hence, \,$\Lambda_{T\,U}$\, is a continuous \,$(\,T,\,U\,)$-controlled $K$-$g$-fusion frame for \,$H$\, with bounds \,$\dfrac{A}{\|\,K\,\|^{\,2}}$\, and \,$B$. \\

$(ii)$\,\,Let \,$\Lambda_{T\,U}$\, be a continuous \,$(\,T,\,U\,)$-controlled $K$-$g$-fusion frame for \,$H$\, with bounds \,$A$\, and \,$B$.\,Since \,$\mathcal{R}\,(\,K\,)$\, is closed, by Theorem \ref{tttt1}, there exists an operator \,$K^{\,\dagger} \,\in\, \mathcal{B}\,(\,H\,)$\, such that \,$K\,K^{\,\dagger}\,f \,=\, f\; \;\forall\; f \,\in\, \mathcal{R}\,(\,K\,)$.\,Then for each \,$f \,\in\, \mathcal{R}\,(\,K\,)$,
\begin{align*}
&\dfrac{A}{\|\,K^{\,\dagger}\,\|^{\,2}}\,\left \|\,f \,\right \|^{\,2} \,\leq\, A\,\|\,K^{\,\ast}\,f\,\|^{\,2}&\\
&\,\leq\, \int\limits_{\,X}\,v^{\,2}\,(\,x\,)\,\left<\,\Lambda_{x}\,P_{F\,(\,x\,)}\,U\,f,\, \Lambda_{x}\,P_{F\,(\,x\,)}\,T\,f\,\right>\,d\mu_{x} \,\leq\, \,B\,\|\,f \,\|^{\,2}. 
\end{align*} 
Thus, \,$\Lambda_{T\,U}$\, is a continuous \,$(\,T,\,U\,)$-controlled \,$g$-fusion frame for \,$\mathcal{R}\,(\,K\,)$\, with bounds \,$\dfrac{A}{\left\|\,K^{\,\dagger}\,\right\|^{\,2}}$\, and $B$. 
\end{proof}

\begin{theorem}
Let \,$K \,\in\, \mathcal{B}\,(\,H\,),\; T,\, U \,\in\, \mathcal{G}\,\mathcal{B}^{\,+}\,(\,H\,)$\, and \,$\Lambda_{T\,U}$\, be a continuous \,$(\,T,\,U\,)$-controlled $K$-$g$-fusion frame for \,$H$\, with frame bounds \,$A,\,B$.\,If \,$V \,\in\, \mathcal{B}\,(\,H\,)$\, with \,$\mathcal{R}\,(\,V\,) \,\subset\, \mathcal{R}\,(\,K\,)$, then \,$\Lambda_{T\,U}$\, is a continuous \,$(\,T,\,U\,)$-controlled $V$-$g$-fusion frame for \,$H$.   
\end{theorem}

\begin{proof}
Since \,$\Lambda_{T\,U}$\, is a continuous \,$(\,T,\,U\,)$-controlled $K$-$g$-fusion frame for \,$H$, for each \,$f \,\in\, H$, we have  
\[A\,\left\|\,K^{\,\ast}\,f\,\right\|^{\,2} \,\leq\, \int\limits_{\,X}\,v^{\,2}\,(\,x\,)\,\left<\,\Lambda_{x}\,P_{\,F\,(\,x\,)}\,U\,f,\, \Lambda_{x}\,P_{\,F\,(\,x\,)}\,T\,f\,\right>\,d\mu_{x} \,\leq\, B\,\|\,f\,\|^{\,2}.\]
Since \,$\mathcal{R}\,(\,V\,) \,\subset\, \mathcal{R}\,(\,K\,)$, by Theorem \ref{th1}, there exists some \,$\lambda \,>\, 0$\, such that \,$V\,V^{\,\ast} \,\leq\, \lambda\, K\,K^{\,\ast}$.\,Thus, for each \,$f \,\in\, H$, we have
\begin{align*}
&\dfrac{A}{\lambda}\,\|\,V^{\,\ast}\,f\,\|^{\,2} \,\leq\, A\,\left\|\,K^{\,\ast}\,f\,\right\|^{\,2}\\
& \,\leq\, \int\limits_{\,X}\,v^{\,2}\,(\,x\,)\,\left<\,\Lambda_{x}\,P_{\,F\,(\,x\,)}\,U\,f,\, \Lambda_{x}\,P_{\,F\,(\,x\,)}\,T\,f\,\right>\,d\mu_{x} \,\leq\, B\,\|\,f\,\|^{\,2}.
\end{align*}
Hence, \,$\Lambda_{T\,U}$\, is a continuous \,$(\,T,\,U\,)$-controlled $V$-$g$-fusion frame for \,$H$.   
\end{proof}

In the following theorem, we will construct a continuous controlled $K$-$g$-fusion frame by using a continuous controlled $g$-fusion frame under some sufficient conditions.

\begin{theorem}\label{th1.1}
Let \,$K \,\in\, \mathcal{B}\,(\,H\,)$\, be an invertible operator on \,$H$\, and \,$\Lambda_{T\,U}$\, be a continuous \,$(\,T,\,U\,)$-controlled $g$-fusion frame for \,$H$\, with frame bounds \,$A,\,B$\, and \,$S_{C}$\, be the associated frame operator.\,Suppose \,$S_{C}^{\,-\, 1}\,K^{\,\ast}$\, commutes with \,$T$\, and \,$U$.\,Then \,$\Gamma_{T\,U} \,=\, \left\{\,\left(\,K\,S_{C}^{\,-\, 1}\,F\,(\,x\,),\, \Lambda_{x}\,P_{F\,(\,x\,)}\,S_{C}^{\,-\, 1}\,K^{\,\ast},\, v\,(\,x\,)\,\right)\,\right\}_{x \,\in\, X}$\, is a continuous \,$(\,T,\,U\,)$-controlled $K$-$g$-fusion frame for \,$H$\, with the corresponding frame operator \,$K\,S_{C}^{\,-\, 1}\,K^{\,\ast}$.  
\end{theorem}

\begin{proof}
Let \,$V \,=\, K\,S_{C}^{\,-\,1}$.\;Then \,$V$\, is invertible on \,$H$\, and \,$V^{\,\ast} \,=\, S_{C}^{\,-\,1}\,K^{\,\ast}$.\,Now, it is easy to verify that 
\begin{equation}\label{eq1.4}
\left\|\,K^{\,\ast}\,f\,\right\|^{\,2} \,\leq\, B^{\,2}\, \left\|\,S_{C}^{\,-\,1}\,K^{\,\ast}\,f\,\right\|^{\,2}\; \;\forall\; f \,\in\, H.
\end{equation}
Now, for each \,$f \,\in\, H$, using Theorem \ref{th1.01}, we have
\begin{align*}
&\int\limits_{\,X}\,v^{\,2}\,(\,x\,)\,\left<\,\Lambda_{x}\,P_{\,F\,(\,x\,)}\,V^{\,\ast}\,P_{\,V\,F\,(\,x\,)}\,U\,f,\, \Lambda_{x}\,P_{F\,(\,x\,)}\,V^{\,\ast}\,P_{V\,F\,(\,x\,)}\,T\,f\,\right>\,d\mu_{x}\\
&=\,\int\limits_{\,X}\,v^{\,2}\,(\,x\,)\,\left<\,\Lambda_{x}\,P_{F\,(\,x\,)}\,V^{\,\ast}\,U\,f,\, \Lambda_{x}\,P_{F\,(\,x\,)}\,V^{\,\ast}\,T\,f\,\right>\,d\mu_{x}\\ 
&=\,\int\limits_{\,X}\,v^{\,2}\,(\,x\,)\,\left<\,\Lambda_{x}\,P_{F\,(\,x\,)}\,U\,S_{C}^{\,-\,1}\,K^{\,\ast}\,f,\, \Lambda_{x}\,P_{F\,(\,x\,)}\,T\,S_{C}^{\,-\,1}\,K^{\,\ast}\,f\,\right>\,d\mu_{x}\\
&\leq\, B\, \|\,S_{C}^{\,-\,1}\,\|^{\,2}\; \left\|\,K^{\,\ast}\,f\,\right\|^{\,2}\\
&\leq\; \dfrac{B}{A^{\,2}}\; \|\,K\,\|^{\,2}\;\|\,f\,\|^{\,2}\; \;[\;\text{using}\; \,B^{\,-1}\,I_{\,H} \,\leq\, S_{C}^{\,-\,1} \,\leq\, A^{\,-1}\,I_{\,H}\;].
\end{align*}
On the other hand, for each \,$f \,\in\, H$, we have
\begin{align*}
&\int\limits_{\,X}\,v^{\,2}\,(\,x\,)\,\left<\,\Lambda_{x}\,P_{\,F\,(\,x\,)}\,V^{\,\ast}\,P_{\,V\,F\,(\,x\,)}\,U\,f,\, \Lambda_{x}\,P_{F\,(\,x\,)}\,V^{\,\ast}\,P_{V\,F\,(\,x\,)}\,T\,f\,\right>\,d\mu_{x}\\ 
&=\,\int\limits_{\,X}\,v^{\,2}\,(\,x\,)\,\left<\,\Lambda_{x}\,P_{F\,(\,x\,)}\,U\,S_{C}^{\,-\,1}\,K^{\,\ast}\,f,\, \Lambda_{x}\,P_{F\,(\,x\,)}\,T\,S_{C}^{\,-\,1}\,K^{\,\ast}\,f\,\right>\,d\mu_{x}\\
&\geq\; A\; \left\|\,S_{C}^{\,-\,1}\,K^{\,\ast}\,f\,\right\|^{\,2} \,\geq\, \dfrac{A}{B^{\,2}}\; \left\|\,K^{\,\ast}\,f\,\right\|^{\,2}\; \;[\;\text{by (\ref{eq1.4})}\;].
\end{align*}
Thus, \,$\Gamma_{T\,U}$\, is a continuous \,$(\,T,\,U\,)$-controlled $K$-$g$-fusion frame for \,$H$.\\Furthermore, for each \,$f \,\in\, H$, we have
\begin{align*}
&\int\limits_{\,X}\,v^{\,2}\,(\,x\,)\,\left<\,\Lambda_{x}\,P_{\,F\,(\,x\,)}\,V^{\,\ast}\,P_{\,V\,F\,(\,x\,)}\,U\,f,\, \Lambda_{x}\,P_{F\,(\,x\,)}\,V^{\,\ast}\,P_{V\,F\,(\,x\,)}\,T\,f\,\right>\,d\mu_{x}\\ 
&=\,\int\limits_{\,X}\,v^{\,2}\,(\,x\,)\,\left<\,\Lambda_{x}\,P_{F\,(\,x\,)}\,U\,S_{C}^{\,-\,1}\,K^{\,\ast}\,f,\, \Lambda_{x}\,P_{F\,(\,x\,)}\,T\,S_{C}^{\,-\,1}\,K^{\,\ast}\,f\,\right>\,d\mu_{x}\\
&=\,\left<\,S_{C}\,S_{C}^{\,-\,1}\,K^{\,\ast}\,f,\, S_{C}^{\,-\,1}\,K^{\,\ast}\,f\,\right> \,=\, \left<\,K\,S_{C}^{\,-\,1}\,K^{\,\ast}\,f,\, f\,\right>.
\end{align*} 
This implies that \,$K\,S_{C}^{\,-\,1}\,K^{\,\ast}$\, is the corresponding frame operator of \,$\Gamma_{T\,U}$.
\end{proof}

In the following theorem, we give a necessary and sufficient condition for continuous controlled $g$-fusion Bessel family to be a continuous controlled \,$K$-$g$-fusion frame with the help of quotient operator.

\begin{theorem}
Let \,$K \,\in\, \mathcal{B}\,(\,H\,)$\, and \,$\Lambda_{T\,U}$\, be a continuous \,$(\,T,\,U\,)$-controlled $g$-fusion Bessel family in \,$H$\, with frame operator \,$S_{C}$.\,Then \,$\Lambda_{T\,U}$\, is a continuous \,$(\,T,\,U\,)$-controlled $K$-$g$-fusion frame for \,$H$\, if and only if the quotient operator \,$\left[\,K^{\,\ast} \,/\, S^{1 \,/\, 2}_{C}\,\right]$\, is bounded.  
\end{theorem}

\begin{proof}
First we suppose that \,$\Lambda_{T\,U}$\, is a continuous \,$(\,T,\,U\,)$-controlled $K$-$g$-fusion frame for \,$H$\, with bounds \,$A$\, and \,$B$.\,Then for each \,$f \,\in\, H$, we have 
\[A\,\|\,K^{\,\ast}\,f\,\|^{\,2} \,\leq\, \int\limits_{\,X}\,v^{\,2}\,(\,x\,)\,\left<\,\Lambda_{x}\,P_{F\,(\,x\,)}\,U\,f,\, \Lambda_{x}\,P_{F\,(\,x\,)}\,T\,f\,\right>\,d\mu_{x} \,\leq\, B\,\|\,f \,\|^{\,2}.\]
Thus, for each \,$f \,\in\, H$, we have
\[A\,\|\,K^{\,\ast}\,f\,\|^{\,2} \,\leq\, \left<\,S_{C}\,f,\, f\,\right> \,=\, \left\|\,S^{1 \,/\, 2}_{C}\,f\,\right\|^{\,2}.\]
Now, it is easy to verify that the quotient operator \,$T \,:\, \mathcal{R}\left(\,S^{1 \,/\, 2}_{C}\,\right) \,\to\, \mathcal{R}\,(\,K^{\,\ast}\,)$\, defined by \,$T\,\left(\,S^{1 \,/\, 2}_{C}\,f\,\right) \,=\, K^{\,\ast}\,f\; \;\forall\; f \,\in\, H$\, is well-defined and bounded.

Conversely, suppose that the quotient operator \,$\left[\,K^{\,\ast} \,/\, S^{1 \,/\, 2}_{C}\,\right]$\, is bounded.\,Then for each \,$f \,\in\, H$, there exists some \,$B \,>\, 0$\, such that
\begin{align*}
&\|\,K^{\,\ast}\,f\,\|^{\,2} \,\leq\, B\,\left\|\,S^{1 \,/\, 2}_{C}\,f\,\right\|^{\,2} \,=\, B\,\left<\,S_{C}\,f,\, f\,\right>\\
&\Rightarrow\, \|\,K^{\,\ast}\,f\,\|^{\,2} \,\leq\, B\,\int\limits_{\,X}\,v^{\,2}\,(\,x\,)\,\left<\,\Lambda_{x}\,P_{F\,(\,x\,)}\,U\,f,\, \Lambda_{x}\,P_{F\,(\,x\,)}\,T\,f\,\right>\,d\mu_{x}.
\end{align*}
Thus, \,$\Lambda_{T\,U}$\, is a continuous \,$(\,T,\,U\,)$-controlled $K$-$g$-fusion frame for \,$H$. 
\end{proof}

Now, we establish that a quotient operator will be bounded if and only if a continuous controlled \,$K$-$g$-fusion frame becomes continuous controlled \,$V\,K$-$g$-fusion frame, for some \,$V \,\in\, \mathcal{B}\,(\,H\,)$. 

\begin{theorem}
Let \,$K \,\in\, \mathcal{B}\,(\,H\,)$\, and \,$\Lambda_{T\,U}$\, be a continuous \,$(\,T,\,U\,)$-controlled $K$-$g$-fusion frame for \,$H$\, with frame operator \,$S_{C}$.\,Let \,$V \,\in\, \mathcal{B}\,(\,H\,)$\, be an invertible operator on \,$H$\, and \,$V^{\,\ast}$\, commutes with \,$T$\, and \,$U$.\,Then the following statements are equivalent:
\begin{description}
\item[$(i)$] \,$\Gamma_{T\,U} \,=\, \left\{\,\left(\,V\,F\,(\,x\,),\, \Lambda_{x}\,P_{F\,(\,x\,)}\,V^{\,\ast},\, v\,(\,x\,)\,\right)\,\right\}_{x \,\in\, X}$\, is a continuous \,$(\,T,\,U\,)$-controlled $V\,K$-$g$-fusion frame for \,$H$. 
\item[$(ii)$]The quotient operator \,$\left[\,\left(\,V\,K\,\right)^{\,\ast} \,/\, S_{C}^{1 \,/\, 2}\,V^{\,\ast} \,\right]$\; is bounded.
\item[$(iii)$]The quotient operator \,$\left[\,\left(\,V\,K\,\right)^{\,\ast} \,/\, \left(\,V\,S_{C}\,V^{\,\ast}\,\right)^{1 \,/\, 2}\,\right]$\, is bounded.
\end{description}
\end{theorem}

\begin{proof}$(i) \,\Rightarrow\, (ii)$\, Suppose \,$\Gamma_{T\,U}$\, is a continuous \,$(\,T,\,U\,)$-controlled \,$V\,K$-$g$-fusion frame with bounds \,$A$\, and \,$B$.\,Then for each \,$f \,\in\, H$, we have 
\begin{align*}
A\, \left\|\,\left(\,V\,K\,\right)^{\,\ast}\,f\,\right\|^{\,2} &\,\leq\, \int\limits_{\,X}\,v^{\,2}\,(\,x\,)\,\left<\,\Lambda_{x}\,P_{\,F\,(\,x\,)}\,V^{\,\ast}\,P_{\,V\,F\,(\,x\,)}\,U\,f,\, \Lambda_{x}\,P_{F\,(\,x\,)}\,V^{\,\ast}\,P_{V\,F\,(\,x\,)}\,T\,f\,\right>\,d\mu_{x}\\
& \,\leq\, B\, \|\,f\,\|^{\,2}.
\end{align*}
By Theorem \ref{th1.01}, for each \,$f \,\in\, H$, we have
\begin{align}
&\int\limits_{\,X}\,v^{\,2}\,(\,x\,)\,\left<\,\Lambda_{x}\,P_{\,F\,(\,x\,)}\,V^{\,\ast}\,P_{\,V\,F\,(\,x\,)}\,U\,f,\, \Lambda_{x}\,P_{F\,(\,x\,)}\,V^{\,\ast}\,P_{V\,F\,(\,x\,)}\,T\,f\,\right>\,d\mu_{x}\nonumber\\
&=\,\int\limits_{\,X}\,v^{\,2}\,(\,x\,)\,\left<\,\Lambda_{x}\,P_{F\,(\,x\,)}\,V^{\,\ast}\,U\,f,\, \Lambda_{x}\,P_{F\,(\,x\,)}\,V^{\,\ast}\,T\,f\,\right>\,d\mu_{x}\nonumber\\ 
&=\,\int\limits_{\,X}\,v^{\,2}\,(\,x\,)\,\left<\,\Lambda_{x}\,P_{F\,(\,x\,)}\,U\,V^{\,\ast}\,f,\, \Lambda_{x}\,P_{F\,(\,x\,)}\,T\,V^{\,\ast}\,f\,\right>\,d\mu_{x} \,=\, \left<\,S_{C}\,V^{\,\ast}\,f,\, V^{\,\ast}\,f\,\right>\label{eq3.6}.
\end{align}
Thus, for each \,$f \,\in\, H$, we have
\[A\, \left\|\,\left(\,V\,K\,\right)^{\,\ast}\,f\,\right\|^{\,2} \,\leq\, \left<\,S_{C}\,V^{\,\ast}\,f,\, V^{\,\ast}\,f\,\right> \,=\, \left\|\,S^{1 \,/\, 2}_{C}\,V^{\,\ast}\,f\,\right\|^{\,2}.\]
Now, we define a operator \,$T \,:\, \mathcal{R}\left(\,S^{1 \,/\, 2}_{C}\,V^{\,\ast}\,\right) \,\to\, \mathcal{R}\left(\,(\,V\,K\,)^{\,\ast}\,\right)$\, by \,$T\,\left(\,S^{1 \,/\, 2}_{C}\,V^{\,\ast}\,f\,\right) \,=\, (\,V\,K\,)^{\,\ast}\,f\; \;\forall\; f \,\in\, H$.\,Then it is easy verify that the quotient operator \,$T$\, is well-defined and bounded. 

$(ii) \,\Rightarrow\, (iii)$\,It is obvious.

$(iii) \,\Rightarrow\, (i)$\, Suppose that the quotient operator \;$\left[\,\left(\,V\,K\,\right)^{\,\ast} \,/\, \left(\,V\,S_{C}\,V^{\,\ast}\,\right)^{1 \,/\, 2}\,\right]$\, is bounded.\,Then for each \,$f \,\in\, H$, there exists \,$B \,>\, 0$\, such that 
\[\left\|\,\left(\,V\,K\,\right)^{\,\ast}\,f\,\right\|^{\,2} \,\leq\, B\, \left\|\,\left(\,V\,S_{C}\,V^{\,\ast}\,\right)^{1 \,/\, 2}\,f\,\right\|^{\,2}.\]
Now, by (\ref{eq3.6}), for each \,$f \,\in\, H$, we have
\begin{align*}
&\int\limits_{\,X}\,v^{\,2}\,(\,x\,)\,\left<\,\Lambda_{x}\,P_{\,F\,(\,x\,)}\,V^{\,\ast}\,P_{\,V\,F\,(\,x\,)}\,U\,f,\, \Lambda_{x}\,P_{F\,(\,x\,)}\,V^{\,\ast}\,P_{V\,F\,(\,x\,)}\,T\,f\,\right>\,d\mu_{x}\\
&\hspace{.45cm}=\,\left<\,S_{C}\,V^{\,\ast}\,f,\, V^{\,\ast}\,f\,\right>  \,=\, \left\|\,\left(\,V\,S_{C}\,V^{\,\ast}\,\right)^{1 \,/\, 2}\,f\,\right\|^{\,2} \,\geq\, \dfrac{1}{B}\,\left\|\,\left(\,V\,K\,\right)^{\,\ast}\,f\,\right\|^{\,2}.
\end{align*}
On the other hand, for each \,$f \,\in\, H$, we have
\begin{align*}
&\int\limits_{\,X}\,v^{\,2}\,(\,x\,)\,\left<\,\Lambda_{x}\,P_{\,F\,(\,x\,)}\,V^{\,\ast}\,P_{\,V\,F\,(\,x\,)}\,U\,f,\, \Lambda_{x}\,P_{F\,(\,x\,)}\,V^{\,\ast}\,P_{V\,F\,(\,x\,)}\,T\,f\,\right>\,d\mu_{x}\\
&=\,\int\limits_{\,X}\,v^{\,2}\,(\,x\,)\,\left<\,\Lambda_{x}\,P_{F\,(\,x\,)}\,U\,V^{\,\ast}\,f,\, \Lambda_{x}\,P_{F\,(\,x\,)}\,T\,V^{\,\ast}\,f\,\right>\,d\mu_{x}\\
&\hspace{.5cm}\leq\, D\,\left\|\,U^{\,\ast}\,f\,\right\|^{\,2} \,\leq\, D\,\|\,U\,\|^{\,2}\,\|\,f\,\|^{\,2}.
\end{align*}
Hence, \,$\Gamma_{T\,U}$\, is a continuous \,$(\,T,\,U\,)$-controlled $V\,K$-$g$-fusion frame for \,$H$.\,This completes the proof. 
\end{proof}

\section{Stability of dual continuous controlled $g$-fusion frame}

\smallskip\hspace{.6 cm}In frame theory, one of the most important problem is the stability of frame under some perturbation.\,P. Casazza and Chirstensen \cite{CC} have been generalized the Paley-Wiener perturbation theorem to perturbation of frame in Hilbert space.\,P. Ghosh and T. K. Samanta \cite{P} discussed stability of dual \,$g$-fusion frame in Hilbert space.\,In this section, we give an important on stability of perturbation of continuous controlled \,$K$-$g$-fusion frame and dual continuous controlled \,$g$-fusion frame.\\

Following theorem provides a sufficient condition on a family \,$\Lambda_{T\,U}$\, to be a continuous controlled \,$K$-$g$-fusion frame, in the presence of another continuous controlled \,$K$-$g$-fusion frame.

\begin{theorem}\label{et1}
Let \,$\Lambda_{T\,U}$\, be a continuous \,$(\,T,\,U\,)$-controlled \,$g$-fusion frame for \,$H$\, and \,$S_{C}$\, be corresponding frame operator.\,Assume that \,$S^{\,-\, 1}_{C}$ is commutes with \,$T$\, and \,$U$.\,Then \,$\Gamma_{T\,U} \,=\, \left\{\,\left(\,S_{C}^{\,-\, 1}\,F\,(\,x\,),\, \Lambda_{x}\,P_{F\,(\,x\,)}\,S_{C}^{\,-\, 1},\, v\,(\,x\,)\,\right)\,\right\}_{x \,\in\, X}$\, is a continuous \,$(\,T,\,U\,)$-controlled $g$-fusion frame for \,$H$\, with the corresponding frame operator \,$S^{\,-\, 1}_{C}$.
\end{theorem}

\begin{proof}
Proof of this theorem directly follows from the Theorem \ref{th1.1}, by putting \,$K \,=\, I_{H}$.
\end{proof}

The family \,$\Gamma_{T\,U}$\, defined in the Theorem \ref{et1} is called the canonical dual continuous \,$(\,T,\,U\,)$-controlled $g$-fusion frame of \,$\Lambda_{T\,U}$.\,We now give the stability of dual continuous controlled $g$-fusion frame.

\begin{theorem}
Let \,$\Lambda_{T\,U}$\, and \,$\Gamma_{T\,U}$\, be two continuous \,$(\,T,\,U\,)$-controlled \,$g$-fusion frames for \,$H$\, with bounds \,$A_{\,1},\, B_{\,1}$\, and \,$A_{\,2},\, B_{\,2}$\, having their corresponding frame operators \,$S_{C}$\, and \,$S_{C^{\,\prime}}$, respectively.\,Consider \,$\Delta_{T\,U} \,=\, \left\{\,\left(\,X\,(\,x\,),\, \Delta_{x},\, v\,(\,x\,)\,\right)\,\right\}_{x \,\in\, X}$\, and \,$\Theta_{T\,U} \,=\, \left\{\,\left(\,Y\,(\,x\,),\, \Theta_{x},\, v\,(\,x\,)\,\right)\,\right\}_{x \,\in\, X}$\, as the canonical dual continuous \,$(\,T,\,U\,)$-controlled $g$-fusion frames of \,$\Lambda_{T\,U}$\, and \,$\Gamma_{T\,U}$, respectively.\,Assume that \,$S^{\,-\, 1}_{C}$\, and \,$S^{\,-\, 1}_{C^{\,\prime}}$\, commutes with both \,$T$\, and \,$U$.\,Then the following statements hold:
\begin{description}
\item[$(i)$]If the condition 
\begin{align*}
&\left|\,\int\limits_{\,X}\,v^{\,2}\,(\,x\,)\,\left(\,\left<\,\Lambda_{x}\,P_{F\,(\,x\,)}\,U\,f,\, \Lambda_{x}\,P_{F\,(\,x\,)}\,T\,f\,\right> \,-\, \left<\,\Gamma_{x}\,P_{G\,(\,x\,)}\,U\,f,\, \Gamma_{x}\,P_{G\,(\,x\,)}\,T\,f\,\right>\,\right)\,d\mu_{x}\,\right|\\
& \,\leq\, D\,\|\,f\,\|^{\,2}
\end{align*}
holds for each \,$f \,\in\, H$\, and for some \,$D \,>\, 0$\, then for all \,$f \,\in\, H$, we have
\begin{align*}
&\left|\,\int\limits_{\,X}\,v^{\,2}\,(\,x\,)\,\left(\,\left<\,\Delta_{x}\,P_{X\,(\,x\,)}\,U\,f,\, \Delta_{x}\,P_{X\,(\,x\,)}\,T\,f\,\right> \,-\, \left<\,\Theta_{x}\,P_{Y\,(\,x\,)}\,U\,f,\, \Theta_{x}\,P_{Y\,(\,x\,)}\,T\,f\,\right>\,\right)\,d\mu_{x}\,\right|\\
& \,\leq\, \dfrac{D}{A_{\,1}\,A_{\,2}}\,\|\,f\,\|^{\,2}.
\end{align*}
\item[$(ii)$]If for each \,$f \,\in\, H$, there exists \,$D \,>\, 0$\, such that
\begin{align*}
&\left|\,\int\limits_{\,X}\,v^{\,2}\,(\,x\,)\,\left<\,T^{\,\ast}\,\left(\,P_{F\,(\,x\,)}\,\Lambda_{x}^{\,\ast}\,\Lambda_{x}\,P_{F\,(\,x\,)} \,-\, P_{G\,(\,x\,)}\,\Gamma_{x}^{\,\ast}\,\Gamma_{x}\,P_{G\,(\,x\,)}\,\right)\,U\,f,\, g\,\right>\,d\mu_{x}\,\right|\\
& \,\leq\, D\,\|\,f\,\|^{\,2}
\end{align*}
then
\begin{align*}
&\left|\,\int\limits_{\,X}\,v^{\,2}\,(\,x\,)\,\left<\,T^{\,\ast}\,\left(\,P_{X\,(\,x\,)}\,\Delta_{x}^{\,\ast}\,\Delta_{x}\,P_{X\,(\,x\,)} \,-\, P_{Y\,(\,x\,)}\,\Theta_{x}^{\,\ast}\,\Theta_{x}\,P_{Y\,(\,x\,)}\,\right)\,U\,f,\, g\,\right>\,d\mu_{x}\,\right|\\
& \,\leq\, \dfrac{D}{A_{\,1}\,A_{\,2}}\,\|\,f\,\|^{\,2}
\end{align*}
\end{description}  
\end{theorem}

\begin{proof}$(i)$
Since \,$S_{C} \,-\, S_{C^{\,\prime}}$\, is self-adjoint so
\begin{align*}
&\left\|\,S_{C} \,-\, S_{C^{\,\prime}}\,\right\| \,=\, \sup\limits_{\|\,f\,\| \,=\, 1}\,\left|\,\left<\,(\,S_{C} \,-\, S_{C^{\,\prime}}\,)\,f,\, f\,\right>\,\right|\\
& \,=\, \sup\limits_{\|\,f\,\| \,=\, 1}\,\left|\,\left<\,S_{C}\,f,\, f\,\right> \,-\, \left<\,S_{C^{\,\prime}}\,f,\, f \,\right>\,\right|\\
&=\, \sup\limits_{\|\,f\,\| \,=\, 1}\,\left|\,\int\limits_{\,X}\,v^{\,2}\,(\,x\,)\,\left(\,\left<\,\Lambda_{x}\,P_{F\,(\,x\,)}\,U\,f,\, \Lambda_{x}\,P_{F\,(\,x\,)}\,T\,f\,\right> \,-\, \left<\,\Gamma_{x}\,P_{G\,(\,x\,)}\,U\,f,\, \Gamma_{x}\,P_{G\,(\,x\,)}\,T\,f\,\right>\,\right)\,d\mu_{x}\,\right|\\
&\leq\, \sup\limits_{\|\,f\,\| \,=\, 1}\,D\,\|\,f\,\|^{\,2} \,=\, D.
\end{align*}
Then 
\begin{align}
\left\|\,S_{C}^{\,-\, 1} \,-\, S_{C^{\,\prime}}^{\,-\, 1}\,\right\| &\,\leq\, \left\|\,S_{C}^{\,-\, 1}\,\right\|\,\left\|\,S_{C} \,-\, S_{C^{\,\prime}}\,\right\|\,\left\|\,S_{C^{\,\prime}}^{\,-\, 1}\,\right\|\nonumber\\
&\leq\, \dfrac{1}{A_{\,1}}\, D\, \dfrac{1}{A_{\,2}} \,=\, \dfrac{D}{A_{\,1}\,A_{\,2}}.\label{eq3.4}
\end{align}
Now, for each \,$f \,\in\, H$, we have
\begin{align*}
&\int\limits_{\,X}\,v^{\,2}\,(\,x\,)\,\left<\,\Delta_{x}\,P_{X\,(\,x\,)}\,U\,f,\, \Delta_{x}\,P_{X\,(\,x\,)}\,T\,f\,\right>\,d\mu_{x}\\
&=\,\int\limits_{\,X}\,v^{\,2}\,(\,x\,)\,\left<\,\Lambda_{x}\,P_{F\,(\,x\,)}\,S_{C}^{\,-\, 1}\,P_{S_{C}^{\,-\, 1}\,F\,(\,x\,)}\,U\,f,\, \Lambda_{x}\,P_{F\,(\,x\,)}\,S_{C}^{\,-\, 1}\,P_{S_{C}^{\,-\, 1}\,F\,(\,x\,)}\,T\,f\,\right>\,d\mu_{x}\\ 
&=\,\int\limits_{\,X}\,v^{\,2}\,(\,x\,)\,\left<\,\Lambda_{x}\,P_{F\,(\,x\,)}\,S_{C}^{\,-\, 1}\,U\,f,\, \Lambda_{x}\,P_{F\,(\,x\,)}\,S_{C}^{\,-\, 1}\,T\,f\,\right>\,d\mu_{x}\\
&=\,\int\limits_{\,X}\,v^{\,2}\,(\,x\,)\,\left<\,\Lambda_{x}\,P_{F\,(\,x\,)}\,U\,S_{C}^{\,-\, 1}\,f,\, \Lambda_{x}\,P_{F\,(\,x\,)}\,T\,S_{C}^{\,-\, 1}\,f\,\right>\,d\mu_{x}\\
&=\,\int\limits_{\,X}\,v^{\,2}\,(\,x\,)\,\left<\,T^{\,\ast}\,P_{F\,(\,x\,)}\,\Lambda^{\,\ast}_{x}\,\Lambda_{x}\,P_{F\,(\,x\,)}\,U\,S_{C}^{\,-\, 1}\,f,\, S_{C}^{\,-\, 1}\,f\,\right>\,d\mu_{x}\\
& \,=\, \left<\,S_{C}\,S^{\,-\, 1}_{C}\,f,\, S^{\,-\, 1}_{C}\,f\,\right> \,=\,  \left<\,f,\, S^{\,-\, 1}_{C}\,f\,\right>.
\end{align*}
Similarly, it can be shown that
\[\int\limits_{\,X}\,v^{\,2}\,(\,x\,)\,\left<\,\Theta_{x}\,P_{Y\,(\,x\,)}\,U\,f,\, \Theta_{x}\,P_{Y\,(\,x\,)}\,T\,f\,\right>\,d\mu_{x} \,=\, \left<\,f,\, S^{\,-\, 1}_{C^{\,\prime}}\,f\,\right>.\]
Therefore, for each \,$f \,\in\, H$, we have
\begin{align*}
&\left|\,\int\limits_{\,X}\,v^{\,2}\,(\,x\,)\,\left(\,\left<\,\Delta_{x}\,P_{X\,(\,x\,)}\,U\,f,\, \Delta_{x}\,P_{X\,(\,x\,)}\,T\,f\,\right> \,-\, \left<\,\Theta_{x}\,P_{Y\,(\,x\,)}\,U\,f,\, \Theta_{x}\,P_{Y\,(\,x\,)}\,T\,f\,\right>\,\right)\,d\mu_{x}\,\right|\\
&=\, \left|\,\left<\,f,\, S^{\,-\, 1}_{C}\,f\,\right> \,-\, \left<\,f,\, S^{\,-\, 1}_{C^{\,\prime}}\,f\,\right>\,\right| \,=\, \left|\,\left<\,f,\, \left(\,S^{\,-\, 1}_{C} \,-\, S^{\,-\, 1}_{C^{\,\prime}}\,\right)\,f\,\right>\,\right|\\
&\,\leq\, \left\|\,S_{C}^{\,-\, 1}\, \,-\, S_{C^{\,\prime}}^{\,-\, 1}\,\right\|\, \|\,f\,\|^{\,2} \,\leq\, \dfrac{D}{A_{\,1}\,A_{\,2}}\,\|\,f\,\|^{\,2}.
\end{align*}
Proof of \,$(ii)$\,\,
In this case, we also find that
\begin{align*}
&\left\|\,S_{C} \,-\, S_{C^{\,\prime}}\,\right\| \,=\, \sup\limits_{\|\,f\,\| \,=\, 1}\,\left|\,\left<\,\left(\,S_{C} \,-\, S_{C^{\,\prime}}\,\right)\,f,\, f\,\right>\,\right|\\
& \,=\, \sup\limits_{\|\,f\,\| \,=\, 1}\,\left|\,\left<\,S_{C}\,f,\, f\,\right> \,-\, \left<\,S_{C^{\,\prime}}\,f,\, f \,\right>\,\right|\\
&=\, \sup\limits_{\|\,f\,\| \,=\, 1}\,\left|\,\int\limits_{\,X}\,v^{\,2}\,(\,x\,)\,\left<\,T^{\,\ast}\,\left(\,P_{F\,(\,x\,)}\,\Lambda_{x}^{\,\ast}\,\Lambda_{x}\,P_{F\,(\,x\,)} \,-\, P_{G\,(\,x\,)}\,\Gamma_{x}^{\,\ast}\,\Gamma_{x}\,P_{G\,(\,x\,)}\,\right)\,U\,f,\, g\,\right>\,d\mu_{x}\,\right|\\
&\leq\, \sup\limits_{\|\,f\,\| \,=\, 1}\,D\,\|\,f\,\|^{\,2} \,=\, D.
\end{align*}
Then for each \,$f \,\in\, H$, we have
\begin{align*}
&\left|\,\int\limits_{\,X}\,v^{\,2}\,(\,x\,)\,\left<\,T^{\,\ast}\,\left(\,P_{X\,(\,x\,)}\,\Delta_{x}^{\,\ast}\,\Delta_{x}\,P_{X\,(\,x\,)} \,-\, P_{Y\,(\,x\,)}\,\Theta_{x}^{\,\ast}\,\Theta_{x}\,P_{Y\,(\,x\,)}\,\right)\,U\,f,\, g\,\right>\,d\mu_{x}\,\right|\\
&=\, \left|\,\left<\,\left(\,S^{\,-\, 1}_{C} \,-\, S^{\,-\, 1}_{C^{\,\prime}}\,\right)\,f,\, f\,\right>\,\right| \,\leq\, \left\|\,S_{C}^{\,-\, 1}\, \,-\, S_{C^{\,\prime}}^{\,-\, 1}\,\right\|\, \|\,f\,\|^{\,2} \,\leq\, \dfrac{D}{A_{\,1}\,A_{\,2}}\, \|\,f\,\|.
\end{align*} 
This completes the proof.
\end{proof}

\end{document}